\documentclass[11pt,a4paper]{article}
\usepackage{amsmath}
\usepackage{amsfonts}
\usepackage{amssymb}
\usepackage{natbib}
\usepackage{amsthm}
\usepackage{array}
\usepackage{graphicx}

\def \R{\mathbb{R}}

\def \N{\mathbb{N}}

\def \E{\mathbb{E}}

\def \Var{\hbox{{\rm Var}}}

\RequirePackage[colorlinks=true,citecolor=blue,linkcolor=blue]{hyperref}

\title{Global optimization using Sobol indices}

\author{Alexandre Janon$^{1,2}$}
\date{1. Laboratoire de Math\'ematiques d'Orsay, Univ. Paris-Sud, CNRS, Universit\'e Paris-Saclay, 91405 Orsay, France.\\2. Inria Saclay - Ile-de-France, B\^at. Turing, Campus de l'Ecole Polytechnique, 91120 Palaiseau, France.\\[.7cm]Homepage: \url{https://www.math.u-psud.fr/\~janon/}\\[.7cm]\today}

\begin{document}
\maketitle
\begin{abstract}
We propose and assess a new global (derivative-free) optimization algorithm, inspired by the LIPO algorithm, which uses variance-based sensitivity analysis (Sobol indices) to reduce the number of calls to the objective function. This method should be efficient to optimize costly functions satisfying the sparsity-of-effects principle.
\end{abstract}

\tableofcontents
\section*{Introduction}
Finding the minimum (or maximum) of a given numerical function (called the objective function) defined over a compact subset of $\R^d$ (the input space) is a fundamental problem with numerous applications such as in model fitting, machine learning or system design.  In general, it requires many evaluations of the objective function over the input space. Each evaluation of the objective can be computationally expensive (for instance, in system design applications, numerical evaluation of the objective can require solving a partial differential equation discretized on a fine mesh), hence the need for  optimization strategies that require as few evaluations as possible. To be efficient, these strategies  rely on assumptions on the objective function and/or input space: convexity or Lipschitz continuity for instance.

In this paper, we will consider an exploration strategy of the input space based on the so-called Sobol global (or variance-based) sensitivity indices. These indices are defined in the following context: the $d$ input parameters are supposed to be independent random variables of known probability distributions, so that the objective function (called the \emph{output function} in this context) is a square-integrable random variable. For each set $u$ of input parameters, we define the (closed) Sobol' index of $u$, which quantifies the fraction of the variance of the output function due to the variability of parameters in $u$. These indices are normalized to be in $[0;1]$, and parameter subsets whose index is close to one have larger "influence" on the output  than subsets with indices close to zero. 

Some output functions may depend on a large number of input parameters, however, computation of the Sobol indices can show that a function verifies a ``sparsity-of-effects'' principles: only a small number of them have a significative effect, or that some interactions of input parameters have small influences on the outputs (for instance, if the output function is a sum of $d$ one-variable functions).

This paper proposes a way to take advantage of this principle to perform optimization. In the first section, we review the definition of Sobol indices and the base sequential optimization strategy; in the second section, we describe the proposed optimization algorithm; and, in the third section, we illustrate the behavior of this algorithm on a numerical test case.

\section{Setup and assumptions}
\subsection{Goal}
We denote by $f: \mathcal D \subset \R^d \rightarrow \R $ the objective function to be minimized. That is, we want to compute:
\[ \min_{\mathcal D} f \]
by using as few as possible evaluations of $f$.

We assume that $\mathcal D = [-1;1]^d$, and we suppose that:
\begin{itemize}
\item $X=(X_1,\ldots,X_d)$ is a uniform random variable on $\mathcal D$,
\item $Y=f(X)$ is a square-integrable real random variable : $Y\in L^2(\mathcal D)$,
\item $Y$ has unit variance: $\Var Y=1$.
\end{itemize}

For $n\in\N^*$, our goal is to define a random sequence $S_n=(X^1, X^2, \ldots, X^n) \in \mathcal D^n$, such that:
\begin{itemize}
\item $S_n$ can be computed by evaluating $f$ exactly $n$ times,
\item for each $i=2,\ldots,n$, $X^{i+1}$ depends on $X^1,\ldots,X^i$ and $f(X^1),\ldots,f(X^i)$,
\item we have that:
\[ \min_{1 \le i \le n} f(X^i) \approx \min_{\mathcal D} f. \]
\end{itemize}

\subsection{Base strategy}
\label{strat}
To compute a ``minimizing'' sequence $S_n=(X^1,\ldots,X^n)$, we take  some subset $\mathcal F$ of functions $\mathcal D \rightarrow \R$ and we use the following algorithm:

\begin{itemize}
\item Inputs: $n\in\N$, $f$ and $\mathcal F$.
\item Initialization : choose $X^1$ uniformly on $\mathcal D$;
\item Iteration: for $i=2,\ldots,n$, repeat:
  \begin{itemize}
    \item choose $X^{i}$ uniformly on:
	 \[ \mathcal D_{i} = \{ x\in\mathcal D \textrm{ s.t. } \exists g \in \mathcal F_i,\; g(x)<\min_{1\le j<i} f(X^j) \} \]
	 where:
	 \[ \mathcal F_i=\{g\in\mathcal F,\; \forall 1\le j<i,\; g(X^j)=f(X^j)\} \]
	\end{itemize}
\end{itemize}

The $\mathcal F$ accounts for ``regularity'' assumptions made on the objective function $f$. At each iteration of the algorithm, the $\mathcal F_i$ set contains functions in $\mathcal F$ that are consistent with the evaluations of $f$ made so far in the algorithm. We look for an input $X^i$ where evaluating $f$ may improve our current minimum, that is, when there exist a function in $\mathcal F_i$ returning a smaller value.

This strategy is used in \citep{malherbe2017global} in the so-called ``LIPO'' algorithm, in the particular case where $\mathcal F$ is the set of the $k$-Lipschitz functions. This paper refers to \citep{hanneke2011rates} and \citep{dasgupta2011two} to call  $\mathcal F_i$  the  active subset of consistent functions.

\subsection{Sobol indices}
For any subset $u\subset\{1,\ldots,d\}$ of size $|u|$,  we denote by $X_u$ the vector:
\[ X_u = (X_k, k\in u)\in[-1;1]^{|u|}. \]
We recall (see for instance the decomposition lemma of \citep{efron1981jackknife}) that for any square integrable function $g:[-1;1]^d\rightarrow\R$ there exists a unique decomposition:
\[ g(X) = \sum_{u\subset\{1,\ldots,d\}} g_u(X_u) \]
where the $g_u(X_u)$ are pairwise orthogonal random variables of $L^2$:
\[ \forall u, v\subset\{1,\ldots,d\}, u\neq v \Longrightarrow \E(g_u(X_u) g_v(X_v))=0. \]

The Sobol index of $g$ relative to the subset of input parameters $u$ is defined by (see \citep{sobol1993,sobol2001global}):
\[ S_u(g) = \frac{\Var(g_u(X_u))}{\Var(g(X))}. \]

Sobol indices can be expressed from the expansion of $g$ on a tensor $L^2$ orthonormal basis: let $(\psi_i)_{i\in\N}$ be an orthonormal basis of $L^2([-1,1])$ (with $\psi_0=\sqrt 2$). Then $g$ can be expanded onto the tensorized basis:
\[ g(X) = \sum_{i\in\N^d, v\subset \{1,\ldots,d\}} a_{i,v}(g) \psi_{i,v}(X_v) \]
where:
\[ \psi_{i,v}(X_v)=\prod_{k=1}^{|v|} \psi_{i_k}(X_{v_k}) \]
for $i=(i_1,\ldots,i_{|v|})\in\N^{|v|}$ and $v=\{v_1,\ldots,v_{|v|}\}$.

Suppose, to simplify, that $\Var(g(X))=1$. Then we have:
\[ S_u(g) = \sum_{i\in\N^{|u|}} a_{i,u}(g)^2 \]
Note that, in the case where the $1D$ basis $(\psi_i)$ is formed by the Legendre polynomials (normalized to have unit variance) we get the expression of the Sobol indices using polynomial chaos expansion (\citep{crestaux2009polynomial}).

\section{Proposed algorithm}
Let's assume that some Sobol indices of $f$ are known to satisfy:
\[ S_u(f) \le s_u,\;\forall u \in \mathcal U \]
where $\mathcal U$ is a set of subsets of $\{1,\ldots,d\}$ and $\{s_u,\;\forall u\in\mathcal U\}$ are known reals in $[0;1]$.

Our minimization algorithm uses the strategy described in \ref{strat} with
\[ \mathcal F = \{ g \in L^2([-1,1]^d) \textrm{ s.t. } \Var g(X)=1,\textrm{ and } \forall u\in\mathcal U, \;S_u(g)\le s_u \} \]

At each iteration, $X^i$ is drawn from uniform distribution on $\mathcal D_i$ using a rejection sampling algorithm : we take a proposal $x$ uniformly sampled on $\mathcal D$, then the existence of $g$ is checked by computing:
\begin{equation}\label{e:defm} m(x)=\min_{h\in\mathcal F_i} h(x) \end{equation}
If 
\[ m(x) < \min_{1\le j<i} f(X^j) \]
then $x$ is accepted and used as $X^i$. Else, $x$ is rejected and a new $x$ is drawn uniformly on $\mathcal D$.

In practice, we use a truncated $L^2$ orthonormal basis of normalized Legendre polynomials (up to degree $D$), and we solve the following optimization problem:
\begin{equation}\label{e:pbmini} m(x)\approx \min_{(a_{k,v})\in\widetilde{\mathcal F_i}} \sum_{k\in\N^d, v\subset\{1,\ldots,d\}} a_{k,v} \psi_{k,v}(x) \end{equation}
with:
\[ \widetilde{\mathcal F_i} = \left\{ (a_{k,v})_{k\in\N^d, v\subset\{1,\ldots,d\}} \textrm{ s.t. }
 \left\{ \begin{array}{l} 
   a_{k,v}=0\; \forall v\subset\{1,\ldots,d\}, \forall k\textrm{ s.t. } \exists \ell,\; k_\ell > D,\\
 	\sum_{k\in\N^d, v\subset\{1,\ldots,d\}} a_{k,v} \psi_{k,v}(X^j) = f(X^j),\; \forall 1\le j<i,\\
	\sum_{k\in\N^d, u\subset\{1,\ldots,d\}, u\neq\emptyset} \sum a_{k,u}^2 \le 1,\\
	\sum_{k\in\N^d} a_{k,u}^2 \le s_u,\; \forall u\in\mathcal U 
	\end{array} \right. \right\} \]
This optimization problem has to be solved multiple times (for various $x$ values), and is affected by the curse of dimensionality; however it is convex (as a quadratically-constrained linear program) and hence can be solved very efficiently without any new evaluation of the $f$ function, e.g.  by an interior point method (see \citep{boyd2004convex} for instance). Hence, this algorithm is computationally interesting if $f$ when sufficiently costly to evaluate. 

\section{Numerical illustrations}

The minimization algorithm above has been implemented in \verb<R<, using the \verb<CVXR< \citep{cvxr2019} package with the \verb<ECOS< solver. We have tested using the $3D$ Rosenbrock function over $[-5;5]^3$:
\[ f(X_1,X_2,X_3) = \frac{1}{26000} \sum_{m=1}^2 100(X_{m+1}-X_m^2)^2 + (1-X_m)^2 \]
(the $\frac{1}{26000}$ being here to ensure that $\Var f\le 1$).

We can notice that no interaction between $X_1$ and $X_3$ occurs.

Instead of using a fixed number $n$ of evaluations of $f$, we set a budget of $N=100$ resolutions of the convex problem \eqref{e:pbmini} and we report the number of evaluations of $f$. We use a maximal polynomial degree of $D=4$. 

We perform four experiments:
\begin{itemize}
\item experiment A : no constraint on the Sobol indices is taken

\item experiment B : we take into account the first-order Sobol indices:
\[ S_{\{1\}} \le 0.42 ,\;\;
   S_{\{2\}} \le 0.46,\;\;
	S_{\{3\}} \le 0.004 \]
as well as the \emph{total} first-order Sobol indices:
\[ S_{\{1\}}+S_{\{1,2\}}+S_{\{1,3\}} \le 0.47,\;
 S_{\{2\}}+S_{\{1,2\}}+S_{\{2,3\}} \le 0.56,\;
 S_{\{3\}}+S_{\{1,3\}}+S_{\{2,3\}} \le 0.06,  \]
Those six indices have been estimated using the \emph{Saltelli estimator} \verb<sobolSalt< of the R package sensitivity \citep{sensi2019}.

\item experiment C : constraints of experiment B, plus zero interaction between variables 1 and 3:
\[ S_{\{1,3\}}=S_{\{1,2,3\}}=0 \]

\item experiment D : only zero interaction between variables 1 and 3 (ie., the constraints of experiment C without those of experiment B).
\end{itemize}

For each experiment, we report:
\begin{itemize}
\item $N_{eval}$ the number of evaluations of $f$ made ;
\item $m$ the computed approximation of $\min f$.
\end{itemize}
The smaller is the better for these two numbers.

Results are given in the following table:

\begin{center}
\begin{tabular}{|c|c|c|}
\hline Experiment & $N_{eval}$ & $m$ \\
\hline A & 93  & 0.0089   \\
\hline B & 78 & 0.0052  \\
\hline C & 44 & 0.0006  \\
\hline D & 45 & 0.0049  \\
\hline
\end{tabular}
\end{center}

We can see that experiments taking Sobol indices into account (B, C and D) improves the results of the algorithm on both criteria ($N_{eval}$ and $m$), and that taking the absence of interaction between $X_1$ and $X_3$ improves further, by halving the necessary number of evaluations to $f$, while yielding a comparable estimate of $\min f$.

\bibliographystyle{plainnat}
\bibliography{biblio,bibi,biblioInt}

\def\cprime{$'$}
\begin{thebibliography}{10}
\providecommand{\natexlab}[1]{#1}
\providecommand{\url}[1]{\texttt{#1}}
\expandafter\ifx\csname urlstyle\endcsname\relax
  \providecommand{\doi}[1]{doi: #1}\else
  \providecommand{\doi}{doi: \begingroup \urlstyle{rm}\Url}\fi

\bibitem[Boyd and Vandenberghe(2004)]{boyd2004convex}
Stephen Boyd and Lieven Vandenberghe.
\newblock \emph{Convex optimization}.
\newblock Cambridge university press, 2004.

\bibitem[Crestaux et~al.(2009)]{crestaux2009polynomial}
T.~Crestaux et~al.
\newblock {Polynomial chaos expansion for sensitivity analysis}.
\newblock \emph{Reliability engineering \& System Safety}, 94\penalty0
  (7):\penalty0 1161--1172, 2009.
\newblock ISSN 0951-8320.

\bibitem[Dasgupta(2011)]{dasgupta2011two}
Sanjoy Dasgupta.
\newblock Two faces of active learning.
\newblock \emph{Theoretical computer science}, 412\penalty0 (19):\penalty0
  1767--1781, 2011.

\bibitem[Efron and Stein(1981)]{efron1981jackknife}
Bradley Efron and Charles Stein.
\newblock The jackknife estimate of variance.
\newblock \emph{The Annals of Statistics}, pages 586--596, 1981.

\bibitem[Fu et~al.(2019)Fu, Narasimhan, Diamond, and Miller]{cvxr2019}
Anqi Fu, Balasubramanian Narasimhan, Steven Diamond, and John Miller.
\newblock \emph{CVXR: Disciplined Convex Optimization}, 2019.
\newblock URL \url{https://CRAN.R-project.org/package=CVXR}.
\newblock R package version 0.99-6.

\bibitem[Hanneke et~al.(2011)]{hanneke2011rates}
Steve Hanneke et~al.
\newblock Rates of convergence in active learning.
\newblock \emph{The Annals of Statistics}, 39\penalty0 (1):\penalty0 333--361,
  2011.

\bibitem[Iooss et~al.(2019)Iooss, Janon, Pujol, with contributions~from
  Baptiste~Broto, Boumhaout, Veiga, Delage, Fruth, Gilquin, Guillaume, {Le
  Gratiet}, Lemaitre, Nelson, Monari, Oomen, Rakovec, Ramos, Roustant, Song,
  Staum, Sueur, Touati, and Weber]{sensi2019}
Bertrand Iooss, Alexandre Janon, Gilles Pujol, with contributions~from
  Baptiste~Broto, Khalid Boumhaout, Sebastien~Da Veiga, Thibault Delage, Jana
  Fruth, Laurent Gilquin, Joseph Guillaume, Loic {Le Gratiet}, Paul Lemaitre,
  Barry~L. Nelson, Filippo Monari, Roelof Oomen, Oldrich Rakovec, Bernardo
  Ramos, Olivier Roustant, Eunhye Song, Jeremy Staum, Roman Sueur, Taieb
  Touati, and Frank Weber.
\newblock \emph{sensitivity: Global Sensitivity Analysis of Model Outputs},
  2019.
\newblock URL \url{https://CRAN.R-project.org/package=sensitivity}.
\newblock R package version 1.16.0.

\bibitem[Malherbe and Vayatis(2017)]{malherbe2017global}
C{\'e}dric Malherbe and Nicolas Vayatis.
\newblock Global optimization of lipschitz functions.
\newblock In \emph{Proceedings of the 34th International Conference on Machine
  Learning-Volume 70}, pages 2314--2323. JMLR. org, 2017.

\bibitem[Sobol(1993)]{sobol1993}
I.~M. Sobol.
\newblock Sensitivity estimates for nonlinear mathematical models.
\newblock \emph{Math. Modeling Comput. Experiment}, 1\penalty0 (4):\penalty0
  407--414 (1995), 1993.
\newblock ISSN 1061-7590.

\bibitem[Sobol(2001)]{sobol2001global}
I.M. Sobol.
\newblock {Global sensitivity indices for nonlinear mathematical models and
  their Monte Carlo estimates}.
\newblock \emph{Mathematics and Computers in Simulation}, 55\penalty0
  (1-3):\penalty0 271--280, 2001.

\end{thebibliography}

\end{document}